\documentclass[11 pt]{article}
\usepackage{fullpage}
\usepackage{amssymb,amsmath,amsfonts, amscd}
\usepackage{array}
\usepackage{amsmath}
\usepackage{amssymb}

\begin{document}

\newtheorem{theorem}{Theorem} [section]
\newtheorem{conjecture}[theorem]{Conjecture}
\newtheorem{tha}{Theorem}
\renewcommand{\thetha}{\Alph{tha}}
\newtheorem{corollary}[theorem]{Corollary}
\newtheorem{lemma}[theorem]{Lemma}
\newtheorem{proposition}[theorem]{Proposition}
\newtheorem{construction}[theorem]{Construction}
\newtheorem{question}[theorem]{Question}
\newtheorem{definition}[theorem]{Definition}
\newtheorem{observation}{Observation}
\newtheorem{remark}[theorem]{Remark}
\newtheorem{fact}[theorem]{Fact}
\setlength{\textwidth}{17cm}
\setlength{\oddsidemargin}{-0.1 in}
\setlength{\evensidemargin}{-0.1 in}
\setlength{\topmargin}{0.0 in}

\def \df {\noindent {\bf Definition. }}

\newcommand{\dist}{{\rm dist}}
\newcommand{\Dist}{{\rm Dist}}
\newcommand{\Forb}{{\rm Forb}}

\newcommand{\qedbox}{$\blacksquare$ \newline}
\newenvironment{proof}%
{%
\noindent{\bf Proof.} } { \hfill\qedbox }

\newcommand{\proofend}{\hfill\qedbox}

\def\qed{\hskip 1.3em\hfill\rule{6pt}{6pt} \vskip 20pt}

\linespread{1.1}
\input epsf
\def\epsfsize#1#2{0.4#1\relax}
\def\O{\text{O}}
\def\o{\text{o}}
\def\ex{\text{ex}}
\def\Z{\mathbb Z}

\def\fl#1{\lfloor #1 \rfloor}
\def\ce#1{\lceil #1 \rceil}

\def \cH{{\cal H }}
\def \cF{{\cal F}}
\def \cG{{\cal G}}
\def \cQ{{\cal Q}}
\def \cA{{\cal A}}
\def \cD{{\cal D}}

\def \f {{\cal F }}
\def \A {{\cal A}}
\def \D {{\cal D}}
\def \fn2 {{\lfloor n/2 \rfloor}}
\def \cn2 {{\lceil  n/2 \rceil}}

\newcommand{\bin}[2]{{#1\choose #2}}
\newcommand{\comp}{\overline}
\newcommand{\exval}{{\rm ex}}

\title{On the strong chromatic number of graphs}
\author{Maria Axenovich\thanks{Department of Mathematics, Iowa
State University, Ames, IA 50011, {\tt
axenovic@math.iastate.edu}}\and Ryan Martin\thanks{Department of
Mathematics, Iowa State University, Ames, IA 50011, {\tt
rymartin@iastate.edu}. Research partially supported by NSA grant
H98230-05-1-0257.}}

\date{}
\maketitle


\begin{abstract}
The strong chromatic number, $\chi_S(G)$,  of an  $n$-vertex graph $G$ is the
smallest number $k$ such that after adding $k\lceil n/k\rceil-n$
isolated vertices to $G$ and considering {\bf any} partition of the
vertices of the resulting graph into disjoint subsets $V_1, \ldots, V_{\lceil
n/k\rceil}$ of size $k$ each, one can find a proper
$k$-vertex-coloring of the graph such that each part $V_i$, $i=1,
\ldots, \lceil n/k\rceil$, contains exactly one vertex of each
color.

For any  graph $G$  with maximum degree $\Delta$, it is easy to
see that $\chi_S(G)\geq\Delta+1$.  Recently, Haxell proved that
 $\chi_S(G) \leq 3\Delta -1$. In this paper, we improve
this bound for graphs with large maximum degree. We show that
$\chi_S(G)\leq 2\Delta$ if $\Delta \geq n/6$ and prove that this
bound is sharp.
\end{abstract}

\section{Introduction}
An $n$-vertex graph $G$ is {\bf strongly $r$-colorable} if
after adding $r\lceil n/r\rceil-n$ isolated vertices to $G$
and considering {\bf any} partition of  the vertices of the resulting graph into disjoint
subsets $V_1, \ldots, V_{\lceil n/r\rceil}$ of size $r$ each, one
can find a proper $r$-vertex-coloring of the  graph such that each
part $V_i$, $i=1, \ldots, \lceil n/r\rceil$, contains exactly one
vertex of each color.
In \cite{F}, it was shown that if a graph
$G$ is strongly $r$-colorable, then it is strongly
$(r+1)$-colorable.

The {\bf strong chromatic number} of $G$,
denoted $\chi_S(G)$, is the smallest positive integer $k$ such
that $G$ is strongly $k$-colorable.

The famous ``cycle plus triangles'' problem of Erd\H{o}s \cite{E},
asking whether $\chi_S(C_{3m})=3$, was answered affirmatively
by Fleischner and Stiebitz \cite{FS}, \cite{FS1},  see also
\cite{S}. In general,  Alon \cite{A}  proved that for any graph
$G$ with maximum degree $\Delta$, $\chi_S(G)\leq c\Delta$, where
$c$ is a very large constant (as the author remarks, $c$ could be reduced to  $10^8$). Recently,  Haxell \cite{H} improved
the bound by Alon drastically, proving that $\chi_S(G)\leq
3\Delta-1$ for any graph $G$ with maximum degree $\Delta$.

As far as the lower bound is concerned, it is easy to see that the
strong chromatic number of a graph with maximum degree $\Delta$ is
at least $\Delta+1$ by taking one of the $V_i$'s to be the
neighborhood of a vertex of maximum degree.

Let $$f(\Delta, n) = \max \{\chi_S(G)\,:\,   G  \mbox{ has
maximum degree }  \Delta  \mbox{ and order } n\}.$$ Therefore, the
best known general bounds are: $$ \Delta +1 \leq f(\Delta, n)\leq
3\Delta -1,$$ for any $\Delta$ and any $n\geq \Delta+1$.

The following theorem is our main result which
gives an exact value for $f(\Delta, n)$ when $\Delta \geq n/6$. It  also provides a
minimum degree condition for the existence of a $K_3$-factor in
tripartite graphs.

\begin{theorem} \label{main}
Let $G$ be a graph on $n$ vertices with maximum degree $\Delta$,
$\Delta \geq n/6$. Then $\chi_S(G)\leq 2\Delta$.
Moreover, for any positive integers $\Delta$ and  $n$,  such that $\Delta \leq n/2$ there is a graph $G_0$ on $n$ vertices,
maximum degree $\Delta$ and $\chi_S(G_0)\geq 2\Delta$.
\end{theorem}

\begin{corollary}
For any positive integer $\Delta$ and any $n$ such that $n/6 \leq  \Delta \leq n/2$,
$f(\Delta, n)=2\Delta$.
Moreover, $f(\Delta, n)\geq 2\Delta$ when  $\Delta\leq n/2$.
\end{corollary}

\section{Proof of Theorem \ref{main}}

In \cite{FS}, \cite{FS1} and others, it was noted  that for specific values of $n$ depending on $\Delta$, there is a graph $G$ such that
$\chi_S(G) \geq 2\Delta$.
We observe here that a similar, but general construction gives the same bound for arbitrary $n$.
Let $\Delta \leq n/2$, let $G_0$ be a graph formed by a disjoint union of
a complete bipartite graph $K_{\Delta, \Delta}$ and $n-2\Delta$ isolated vertices.
Assume that $\chi_S(G_0)\leq 2\Delta-1$, that is,  for $r=2\Delta -1$, any partition of  $V(G_0)$ and $r\lceil n/r\rceil-n$ isolated vertices  into
$t=\lceil n/r\rceil$ sets of equal sizes, $V_1, \ldots, V_t$,
allows a proper $r$-coloring of the resulting graph such that each $V_i$ uses all the colors.
Note that $t \geq \lceil 2\Delta/(2\Delta-1)\rceil =2$.
Now, let $A, B$ be the partite sets of a complete bipartite subgraph of $G_0$ with $|A|=|B|=\Delta$,
let $A\subseteq V_1$ and $B\subseteq V_2$. Then it is easy to see that it is impossible to find the
desired $r$-coloring.

Together with the upper bound which we prove below,
we shall have that $\chi_S(G_0)=2\Delta$ when $ n/6 \leq \Delta \leq n/2$. \\

Now we shall prove the main statement of Theorem \ref{main} by
providing an upper bound on the strong chromatic number.
Let $G$ be a graph on $n$ vertices with maximum degree $\Delta\geq n/6$.

Let $ \Delta \geq n/2$. Then $2\Delta \geq n$ and
we trivially have that $ \chi_S(G) \leq n \leq 2\Delta$.

Let  $n/4\leq \Delta< n/2$.  Thus, $n/2\leq 2\Delta<n$ and
we have to partition $V(G)$ and needed isolated vertices arbitrarily into
two sets $V_1$ and $V_2$, $|V_1|=|V_2|$.
Each vertex in $V_1$ is nonadjacent to at least $|V_2|/2$ vertices in
$V_2$ and vice versa.  Consider the bipartite complement $G'$ of
this graph.  That is, the edge set of $G'$ consists of all  pairs
$\{v_1,v_2\}$, $v_1\in V_1$ and $v_2\in V_2$ such that $\{v_1,
v_2\}\notin E(G)$.
We claim that for each $S\subseteq V_1$, $|N(S)|\geq |S|$.
Indeed, assume that there is a set $S'\subseteq V_1$ for which $|N(S')|<|S'|$.
 We have then that $|S'|>|V_1|/2$, thus, for any vertex
$v\in V_2\setminus N(S')$, $v$ is adjacent to at most $|V_1|-|S'|< |V_1|/2$
vertices, a contradiction.
Applying the K\"onig-Hall theorem to $G'$ gives
a perfect matching, which provides a proper coloring of the
original graph, $G$, with $2\Delta$  colors, each represented
exactly once in $V_1$ and exactly once in $V_2$.

Let  $n/6\leq  \Delta < n/4$. As before, in order to verify that
 $\chi_S(G)\leq r=2\Delta$,
we need to add  $r\lceil n/r\rceil-n$ isolated
vertices to $G$ and partition the resulting vertex set arbitrarily
into  parts  $V_1, V_2, V_3$ of equal sizes.
We shall be treating this case by analyzing and  extending partial colorings.

A {\bf partial strong coloring} of $G$ with respect to $V_1, V_2, V_3$  is a
proper coloring of a subset of the vertices of $G$ such that no
two colored vertices in the same part $V_i$, $i=1, 2, 3$  have the same color and each
color class contains exactly $3$ vertices. For a set $S$  of
vertices and a vertex coloring $\chi$, we say that $S$ is {\bf
partially multicolored} by $\chi$ if any two vertices in $S$,
which are colored by $\chi$, have distinct colors.
Let $\chi$ be a maximal partial strong coloring of $G$ with respect to $V_1, V_2, V_3$.
We will show that we can always {\bf enlarge} such partial strong coloring; i.e., create
another partial strong coloring with more colors, until we color
all the vertices.  For a color $c$, we denote the vertices of this
color $\{c_1,c_2,c_3\}$, where $c_i\in V_i$ for $i=1,2,3$. We fix
$v_1\in V_1, v_2\in V_2, v_3\in V_3$ such that none of $v_1, v_2,
v_3$ are colored by $\chi$.  For $i=1,2,3$, define the following
set:
$$ X_i \stackrel{\rm def}{=}
   \{u \in V_i :  v_i \mbox{ is not adjacent to a
                             vertex of color } \chi(u)\} \cup
         \{u\in V_i: u \mbox{ is not colored by } \chi\} . $$

Observe that any colored vertex in $X_i$ can be replaced by $v_i$, $i=1,2,3$ to
create another strong partial coloring.
Note also that  $$|X_i|\geq |V_i| - \deg(v_i)+t_i \geq \Delta, \quad i=1,2,3,$$
where $t_i$ is the number of neighbors of $v_i$ in $\{v_1, v_2, v_3\}$.

To simplify the notation, we shall assume that no color of $\chi$
is labelled by  $x,v$, or $w$, we reserve $x_i$ or $w_i$ to denote
a vertex in $X_i$ (it might be colored or not colored), and $v_i$
are the vertices fixed above. We shall write $z\sim y$, $z\not\sim
y$ if $zy\in E(G)$, $zy \notin E(G)$, respectively. For disjoint subsets
$S_1$, $S_2$ of vertices of $G$ and a vertex $z$, $z\notin S_1$,  we write
$S_1\sim S_2$ if each vertex in $S_1$ is adjacent to all vertices
in $S_2$, $S_1 \not\sim S_2$ if there are no edges between $S_1$
and $S_2$, and $z\sim S_1$, $z\not\sim S_1$ if $\{z\}\sim S_1$,
$\{z\}\not\sim S_1$, respectively.

To start the proof, we give two  lemmas which allow us either to enlarge $\chi$
 or to replace $\chi$  with another
partial strong coloring such that some three specific vertices
become uncolored and the number of colors remains the same.

\renewcommand{\theenumi}{(\arabic{enumi})}
\begin{lemma}\label{lem:uncolor}
Let $x_i\in X_i$, $i=1,2,3$. If  $\{x_1, x_2, x_3\}$ is partially multicolored
then there is a strong partial coloring with as many color classes as
$\chi$ and with $x_i$'s being uncolored.
\end{lemma}

\begin{proof}
Suppose each $x_i$, $i=1,2,3$ is colored; i.e.,  $x_1=a_1$,
$x_2=b_2$, $x_3=c_3$ with distinct colors $a, b, c$. Replace color
classes $a$, $b$ and $c$ with new color classes $\{v_1, a_2,
a_3\}$, $\{b_1, v_2, b_3\}$ and $\{c_1, c_2, v_3\}$, see Figure \ref{fig:lem211}.

\begin{figure}[ht]
\begin{center}
\epsffile{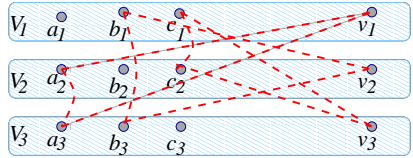}
\end{center}
\caption{Color switches for Lemma \ref{lem:uncolor}}
\label{fig:lem211}
\end{figure}

Now, suppose exactly one $x_i$ is uncolored,  without loss of
generality, $x_1=a_1$, $x_2=b_2$ and $x_3$ is uncolored.  Replace
color classes $a$ and $b$ with new color classes $\{v_1, a_2,
a_3\}$, $\{b_1, v_2, b_3\}$.
Suppose exactly one $x_i$ is colored, without loss of
generality, $x_1=a_1$ and $x_2, x_3$ are uncolored.  Replace color
class $a$ with the new color class $\{v_1, a_2, a_3\}$.
Each case makes $\{x_1, x_2, x_3\}$ uncolored.
\end{proof}

\begin{lemma} $\;$ \label{lem:extend}
\begin{enumerate}
\item If there is a set $\{x_1,x_2,x_3\}$, with
   $x_i\in X_i$, $i=1,2,3$, which induces an independent set
   and is partially multicolored  then $\chi$ can be enlarged.
   \label{lem:extend1}
\item If there is a set $\{x_i,x_i',x_j,x_k\}$,
   with $x_i,x_i'\in X_i$, $x_j\in X_j$, $x_k\in X_k$,
   $\{i,j,k\}=\{1,2,3\}$ such that $\{x_j,x_k\}$ is
   partially multicolored, and both $\{x_i,x_j,x_k\}$ and
   $\{x_i',x_j,x_k\}$ induce independent sets, then $\chi$ can be
   enlarged. \label{lem:extend2}
\item Let a set $\{x_1,x_2,x_3\}$, with $x_i\in X_i$, $i=1,2,3$,
   induce an independent set and the set $\{v_1, v_2, v_3\}$ induce
   neither an independent set nor a clique. Then either  $\chi$ can
   be enlarged or one can find another partial strong coloring with
   as many color classes as in $\chi$ and with three uncolored
   vertices $x'_i\in X_i$, $i=1, 2, 3$ that induce a $K_3$.
   \label{lem:extend3}
\end{enumerate}
\end{lemma}

\begin{proof}
\noindent \underline{\it \ref{lem:extend1}}: By
Lemma~\ref{lem:uncolor} there is a partial
strong coloring with as many color classes as in $\chi$ and such
that $x_1,x_2,x_3$ are uncolored. We can give these vertices a
new color thus enlarging the coloring. \\

\noindent \underline{\it \ref{lem:extend2}}: If either
$\{x_i,x_j,x_k\}$ or $\{x_i',x_j,x_k\}$ is partially multicolored
then we can use \ref{lem:extend1}, otherwise assume, without loss
of generality, that $i=1,j=2,k=3$ and $x_1=a_1$, $x_1'=b_1$,
$x_2=b_2$, $x_3=a_3$, for distinct colors $a, b$. Consider the
following sets of vertices: $\{v_1, b_2, a_3\}$, $\{b_1, v_2,
b_3\}$ and $\{a_1, a_2, v_3\}$. They are independent because of
the  definition of $X_i$'s, $i=1,2,3$.  We can color vertices in
each of these sets with the same new color, which replaces color
classes $a$, $b$ and saturates  vertices $\{v_1,v_2,v_3\}$, thus
enlarging $\chi$.  \\

\noindent \underline{\it \ref{lem:extend3}}:  We can assume,
without loss of generality, that
$v_1\sim v_2$ and $v_2\not\sim v_3$. \\

{\it Case 1.} $\chi(x_1)= \chi(x_2)=\chi(x_3)=a.$
Replace the color class $a$ with two new color classes:  $\{x_1, v_2, v_3\}$ and $\{v_1, x_2, x_3\}$,
thus enlarging $\chi$.\\

{\it Case 2.} $\chi(x_2)=\chi(x_3)=a$.
If $x_1$ is not colored by $\chi$,  replace color class $a$ with two new color classes: $\{a_1, v_2, v_3\}$ and $\{x_1, x_2, x_3\}$.
If $x_1$ is colored $b$, replace color classes $a$ and $b$ with the following three new color classes:
$\{a_1, v_2, v_3\}$,  $\{v_1, b_2, b_3\}$,  $\{x_1, x_2, x_3\}$.
This enlarges $\chi$.\\

{\it Case 3.} $\chi(x_1)=\chi(x_2)=a$.
If $x_3\sim \{v_1, v_2\}$ then replace $x_3$ with $v_3$ in its color class
if $x_3$ is colored by $\chi$. Then $v_1, v_2, x_3$ are three uncolored
vertices inducing a clique. Let $\{x_1', x_2', x_3'\}= \{v_1, v_2, x_3\}$.
If $x_3\not\sim v_1$  then replace $x_3$ by $v_3$ in its color class (if $x_3$ is colored) and
replace a color class $a$ with two new color classes: $\{v_1, x_2, x_3\}$, $\{x_1, v_2, a_3\}$,
thus enlarging $\chi$.
If $x_3\not\sim v_2$  then replace $x_1$ by $v_1$ in its color class,
replace $x_3$ by $v_3$ in its color class (if $x_3$ is colored), and give a new color
to the independent set $\{x_1, v_2, x_3\}$, thus enlarging $\chi$.
Note that the case when $\chi(x_1)=\chi(x_3)=a$ is symmetric.\\

{\it Case 4.} $\{x_1,x_2, x_3\}$ is partially multicolored.
This is  part \ref{lem:extend1} of this Lemma.

\end{proof}

\begin{figure}[ht]
\begin{center}
\epsffile{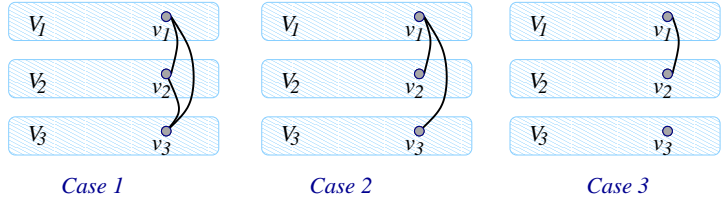}
\end{center}
\caption{Cases 1, 2, 3} \label{fig:cases}
\end{figure}

Next, we consider three cases depending on how many edges the set
$\{v_1,v_2,v_3\}$ induces in $G$.  We shall greedily choose
appropriate $x_i\in X_i$, $i=1,2,3$ and enlarge the coloring.
The proof begins with Case 1, where $\{v_1,v_2,v_3\}$ induces three
edges.  In this case, the coloring can be enlarged.  In Case 2,
$\{v_1,v_2,v_3\}$ induces two edges, without loss of generality
$v_2\not\sim v_3$, and either the coloring can be enlarged or
another coloring with the same number of colors can be found so
that there are three pairwise adjacent uncolored vertices reducing
the analysis to  Case 1.  Finally, in Case 3 there is only one
edge, without loss of generality $v_1v_2$,  induced by
$\{v_1,v_2,v_3\}$.  In this case, either the coloring can be
enlarged or we can find a coloring that puts us in Case 2 or Case 1.
See Figure \ref{fig:cases}.

In Cases 1 and 2 we shall need the following parameter:
$$ q\stackrel{\rm def}{=}
   \max \{|N(x)\cap X_j|: \quad  x\in X_i;~
          i\neq j\mbox{ with }i,j \in \{1,2,3\}\}. $$

\centerline{\large\bf Case 1: $v_1\sim v_2$, $v_1\sim v_3$ and
$v_2\sim v_3$} \vspace{.125in}

We have $|X_i|\geq |V_i|-(\deg(v_i)-2)\geq\Delta+2$, for
$i=1,2,3$. Without loss of generality, assume that $q=|N(x_1)\cap
X_2|$, for $x_1\in X_1$.  Let  $x_2\in X_2\setminus N(x_1)$, be a
vertex not of color $\chi(x_1)$. Consider
$S=X_3\setminus\left(N(x_1)\cup N(x_2)\right)$.  By the choice of
$x_1$, $|S|\geq |X_3|-(\Delta-q)-q\geq (\Delta+2)-\Delta=2$, thus
there are two vertices $x_3,x_3'\in X_3$ nonadjacent to both $x_1$
and $x_2$. Therefore, Lemma~\ref{lem:extend}~\ref{lem:extend2} can
be applied to the four vertices $x_1,x_2,x_3,x_3'$ to enlarge the
coloring.
\\~\\

\centerline{\large\bf Case 2: $v_1\sim v_2$, $v_1\sim v_3$ and
$v_2\not\sim v_3$} \vspace{.125in}

In this case, $|X_1|\geq\Delta+2$ and $|X_2|,|X_3|\geq\Delta+1$.
Let $q= |N(x_i)\cap X_j|$, $x_i\in X_i$ and
let $k\in \{1,2,3\}\setminus\{i,j\}$.
Let $x_j\in X_j\setminus N(x_i)$, let $x_k\in X_k \setminus
(N(x_i)\cup N(x_j)).$
Note that such $x_j$ and $x_k$ exist since $|X_j|\geq \Delta+1$ and
$|X_k\ \setminus (N(x_i)\cup N(x_j))|\geq \Delta+1 - q - (\Delta -q) \geq 1$.

Therefore, we can apply Lemma \ref{lem:extend}~\ref{lem:extend3} to an independent set $\{x_i, x_j, x_k\}$. This either enlarges $\chi$ or reduces Case 2 to Case 1.\\

\centerline{\large \bf Case 3: $v_1\sim v_2$, $v_1\not\sim v_3$ and
$v_2\not\sim v_3$} \vspace{.125in}

\begin{figure}[ht]
\begin{center}
\epsffile{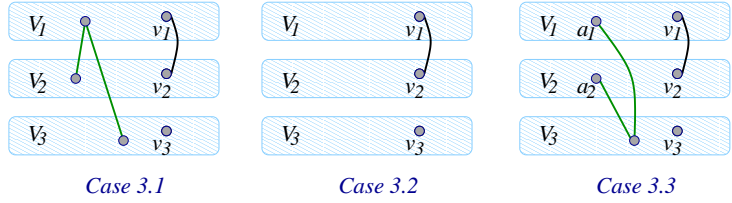}
\end{center}
\caption{Subcases of Case 3}\label{fig:cases3}
\end{figure}

We show that in each of the Cases 3.1--3.3 one can enlarge the
coloring, either directly or by finding a coloring with the same
number of colors that satisfies either the conditions of Case 2 or
the conditions of Case 1. These subcases are arranged according  to the
presence of specific paths in $X_1\cup X_2\cup X_3$,  see Figure \ref{fig:cases3}.

~\\

\noindent {\bf  Case 3.1.} There is no path with three vertices
$w_1,w_2,w_3$;  $w_i\in X_i$, $i=1,2,3$.  \\

We have that  $|X_1|, |X_2|\geq \Delta+1$ and
$|X_3|\geq \Delta$. Let $G_{i,j}$ be the bipartite subgraph of
$G$ induced by the edges of $G$
between $X_i$ and $X_j$, with $i\neq j$ and $i,j\in\{1,2,3\}$.
Note that $G_{i,j}=G_{j,i}$. Moreover, the distinct graphs $G_{i,j}$
are pairwise vertex-disjoint.
If one of $G_{i,j}$ has a nonedge $x_i\not\sim x_j$, then for any $x_k\in X_k$, $k\in \{1,2,3\}\setminus \{x_i, x_j\}$,
$\{x_1, x_2, x_3\}$ is an independent set.
Thus, we can assume that each $G_{i,j}$ is a complete bipartite graph. It is easy to see that in this case there is also an independent set $\{x_1, x_2, x_3\}$,
$x_i\in X_i$, $i=1, 2, 3$. Now,  we can apply Lemma~\ref{lem:extend}~\ref{lem:extend3} to  $\{x_1, x_2, x_3\}$ and either enlarge the coloring or reduce the analysis to Case 1.\\

\noindent {\bf Case 3.2.} There is a path $P$ with three vertices
$w_1,w_2,w_3$;  $w_i\in X_i$, $i=1,2,3$ such that either the
vertices of $P$ are  partially multicolored or  the middle vertex
of $P$ is in $X_1\cup X_2$. \\

If $P$ is partially multicolored, we can apply
Lemma~\ref{lem:uncolor} immediately to obtain a
partial strong coloring with as many colors as $\chi$ and with
vertices of $P$ being uncolored.  We can now choose $v_i=w_i$,
$i=1,2,3$ and use Case 2 or Case 1.

If $P$ has repeated colors on its
vertices, these can be only endvertices of $P$.  Without loss
of generality, let the midpoint of $P$ be $w_1\in X_1$, let
$a_2=w_2$, $a_3=w_3$ be the endpoints of $P$.
If $w_1$ is not colored, replace color class $a$ with an independent set $\{a_1, v_2, v_3\}$.
If $w_1$ has color $b$, then, in addition, replace the color class $b$ with an
independent set $\{v_1, b_2, b_3\}$.
This uncolors $w_1, w_2, w_3$ and brings us to Case 2. \\

\noindent {\bf Case 3.3.} There is a path $(w_1,w_3,w_2)$; with
$w_i\in X_i$, for $i=1,2,3$, and $w_1,w_2$ of the same color.
Moreover, there are no paths satisfying the conditions of Case 3.2.\\

Note that there is no independent set $\{x_1, x_2, x_3\}$, $x_i
\in X_i$, otherwise we can either enlarge the coloring or reduce
the analysis to Case 1 by Lemma
\ref{lem:extend}~\ref{lem:extend3}. Note also that if $x_1\sim x_2$,
$x_i\in X_i$, $i=1,2$, then $\{x_1, x_2\}\not\sim X_3$, otherwise it is Case 3.2. Therefore,
we have that the bipartite subgraph of $G$ with parts $X_1, X_2$
induces  one nontrivial connected component $F$ which  must be a
complete bipartite graph. Since $v_i\in X_i$, $i=1,2,3$, and
$v_1\sim v_2$, $v_1, v_2\in V(F)$. Let $B_1\subseteq X_1,
B_2\subseteq X_2$, be the partite sets of $F$. Let
$A_i=X_i\setminus B_i$, $i=1,2$. Then we have that $B_1\sim B_2$,
$A_1\cup A_2\sim X_3$, $A_1\not\sim A_2$. Then, in particular, we
have that $a_i\in A_i$, $i=1,2$, and $|A_1|=|A_2|=1$, otherwise we
shall find a path satisfying Case 3.2. Since $|X_1|, |X_2|\geq
\Delta+1$, we have that $|B_1|=|B_2|=\Delta$.
Therefore, we can conclude that $|X_1|=|X_2|=\Delta+1$ and $|X_3|=\Delta$.\\

\noindent {\it Claim.}  The vertices $v_1,v_2,v_3$ are the
only uncolored vertices and every color class other than $a$ has
exactly one member in  $X_1\cup X_2$.
\begin{quote}
   {\it Proof of Claim}
   Let $b$ be a color used by $\chi$, $b\neq a$, not present on vertices of $X_1$.
   $N(v_2)=B_1$, so $v_2\not\sim
   b_1$ and $v_2\not\sim b_3$. This  implies that $b_2\in X_2$.
   Thus, any color $b$, $b\neq a$, is used on some vertex in
   $X_1\cup X_2$.

   Let $t$ be the number of uncolored vertices in each $V_i$,
   $i=1,2,3$; i.e., the number of color classes in $\chi$ is
   $2\Delta-t$.  The fact that each color class other than $a$
   contains at least one member of $X_1\cup X_2$ and $a$
   contains two such members gives  that  $|X_1|+|X_2|\geq
   (2\Delta-t+1)+2t$.  Here, the expression in parenthesis
   gives the lower bound on number of colored vertices in $X_1$
   and $X_2$ and $2t$ is the number of uncolored vertices in
   $X_1$ and $X_2$.  Because $|X_1|+|X_2|=2\Delta+2$, we have
   that $t=1$.  As a result, every vertex other than
   $v_1,v_2,v_3$ is colored and every color class other than
   $a$ contains exactly one vertex from $X_1\cup X_2$.
\end{quote}

By Claim, there are $2\Delta-2$ colors different from $a$ in
$\chi$.  Let $\nu$ be the number of neighbors of $v_3$ colored
differently than $a$.  Since
$v_3\sim\{a_1,a_2\}$,  we have that
$\deg(v_3)\geq\nu+2$.  For a color $c$, $c\neq a$, the conditions
$v_3\not\sim c_1$ and $v_3\not\sim c_2$  imply that $c_3\in X_3$.
Using also the fact that $v_3\in X_3$, we have that
$ |X_3|\geq (2\Delta-2-\nu)+1. $  Since $|X_3|=\Delta$, we have that
$(2\Delta-2-\nu)+1\leq\Delta$, thus $\nu\geq\Delta-1$. Therefore,
$\deg(v_3)\geq\nu+2\geq\Delta+1$, a contradiction.

This concludes Case 3.3, and the proof of Theorem \ref{main}. \proofend

\section{Concluding Remarks}

It should be noted that Theorem \ref{main} is equivalent to the
following:

\begin{corollary}
Let $G$ be a tripartite graph with parts of size $n$ each.
If the minimum degree of $G$ is at least $3n/2$ then
$G$ has a $K_3$-factor.
\end{corollary}

This result provides another sufficient condition for the existence
of $K_3$-factors. For other results in this area,  see for
example, \cite{CH,HSz,AY,Fi,MM,MSz, J}.
It also came to author's attention after the manuscript has been
submitted that this problem has been considered independently
in \cite{JJM}, by  treating  $r$-factors in multipartite graphs under maximum degree conditions.

\section{Acknowledgments}
The authors are indebted to anonymous referees for their
careful work. They are especially thankful for
a referee, who provided Lemma \ref{lem:extend}~\ref{lem:extend3},
which helped shortening the proofs.


\begin{thebibliography}{99}

\bibitem{A}  N. Alon, The strong chromatic number of a graph.
{\it Random Structures Algorithms} {\bf 3} (1992), no. 1, 1-7.

\bibitem{AY}  N. Alon and R. Yuster, $H$-factors in dense graphs.
{\it J. Combin. Theory Ser. B} {\bf 66} (1996), no. 2, 269-282.

\bibitem{CH}  K. Corr\'adi and A. Hajnal, On the maximal number
of independent circuits in a graph.  {\it Acta Math. Acad. Sci.
Hungar.} {\bf 14} (1963) 423-439.

\bibitem{E} P. Erd\H{o}s, On some of my favourite problems in
graph theory and block designs.  Graphs, designs and combinatorial
geometries (Catania, 1989). {\it Matematiche (Catania)} {\bf 45}
(1990), no. 1, 61-73 (1991).

\bibitem{F} M. Fellows, Transversals of vertex partitions in
graphs.  {\it SIAM J. Discrete Math.} {\bf 3} (1990), no. 2,
206-215.

\bibitem{Fi} E. Fischer, Variants of the Hajnal-Szemer\'edi
theorem.  {\it J. Graph Theory} {\bf 31} (1999), no. 4, 275-282.

\bibitem{FS} H. Fleischner and M. Stiebitz, A solution to a
colouring problem of P. Erd\H{o}s.  Special volume to mark the
centennial of Julius Petersen's ``Die Theorie der regulren
Graphs'', Part II.  {\it Discrete Math.} {\bf 101} (1992), no.
1-3, 39-48.

\bibitem{FS1} H. Fleischner and M. Stiebitz, Some remarks on the
cycle plus triangles problem. {\it The mathematics of Paul
Erd\H{o}s, II},  136-142, Algorithms Combin., 14, Springer,
Berlin, 1997.


\bibitem{J} R. Johansson, Triangle factors in a balanced blown-up triangle.
{\it Discrete Mathematics}, {\bf 211} (2000), 249-254.



\bibitem{JJM} A. Johansson,  R. Johansson, K. Markstr\"om,
Factors of $r$-partite graphs, {\it personal communication}.

\bibitem{HSz} A. Hajnal and E. Szemer\'edi, Proof of a conjecture
of P. Erd\H{o}s.  {\it Combinatorial theory and its applications,
II (Proc. Colloq., Balatonf\"{u}red, 1969)}, pp. 601-623.  {\it
North-Holland, Amsterdam}, 1970.

\bibitem{H} P. E. Haxell, On the strong chromatic number,
{\it Combin. Probab. Comput.} {\bf 13} (2004), no. 6, 857-865.

\bibitem{MM} Cs. Magyar and R. Martin, Tripartite version of the
Corr\'{a}di-Hajnal theorem.  {\it Discrete Math.} {\bf 254}
(2002), no. 1-3, 289-308.

\bibitem{MSz} R. Martin and E. Szemer\'{e}di, Quadripartite
version of the Hajnal-Szemer\'edi Theorem, submitted.

\bibitem {S} H. Sachs,  Elementary proof of the
cycle-plus-triangles theorem.  Combinatorics, Paul Erd\H{o}s is
eighty, Vol. 1,  347-359, Bolyai Soc. Math. Stud., J\'anos Bolyai
Math. Soc., Budapest, 1993.

\end{thebibliography}
\end{document}